\documentclass{article}
\usepackage{amsmath,amssymb}
\usepackage{amsthm}
\usepackage{enumerate}
\makeatletter
\@addtoreset{equation}{section}
\makeatother
\usepackage[text={28cc,42cc},centering]{geometry}

\newcommand\cl[1][G]{\ensuremath{\mathrm{cl}(#1)}}
\newcommand\Cf[1][i]{\ensuremath{C_5^{(#1)}}}
\newcommand\VC[1][i]{\ensuremath{V(C_5^{(#1)})}}
\newcommand\se{\subseteq}
\newcommand\beq[1]{\begin{equation}\label{#1}}
\newcommand\eeq{\end{equation}}
\newtheorem*{MT}{Main Theorem}
\newtheorem{theorem}{Theorem}[section]
\newtheorem{lem}{Lemma}[section]
\newtheorem{cor}{Corollary}[section]
\theoremstyle{definition}
\newtheorem{definition}{Definition}[section]

\begin{document}

\author{Nikolay Kolev, Nedyalko Nenov}
\title{AN EXAMPLE OF A 16-VERTEX FOLKMAN EDGE (3,4)-GRAPH WITHOUT 8-CLIQUES}

\maketitle

\begin{abstract}
In \cite{NN:6} we computed the edge Folkman number $F(3,4;8)=16$. There we used and
announced without proof that in any blue-red coloring of the edges of the graph
$K_1+C_5+C_5+C_5$ there is either a blue 3-clique or red 4-clique. In this
paper we give a detailed proof of this fact.


\textbf{2000 Math.\ Subject Classification.} 05C55
\end{abstract}

\section{Introduction}

Only finite non-oriented graphs without multiple edges and loops are considered.
We call a $p$-clique of the graph $G$ a set of $p$ vertices each two of which
are adjacent. The largest positive integer $p$ such that $G$ contains
a $p$-clique is denoted by \cl. A set of vertices of the graph $G$
none two of which are adjacent is called an independent set. In this paper
we shall also use the following notations:
\begin{itemize}
  \item $V(G)$ is the vertex set of the graph $G$;
  \item $E(G)$ is the edge set of the graph $G$;
  \item $N(v)$, $v\in V(G)$ is the set of all vertices of $G$ adjacent to $v$;
  \item $G[V]$, $V\se V(G)$ is the subgraph of $G$ induced by $V$;
  \item $\chi(G)$ is the chromatic number of $G$;
  \item $K_n$ is the complete graph on $n$ vertices;
  \item $C_n$ is the simple cycle on $v$ vertices.
\end{itemize}

The equality $C_n=v_1v_2\dots v_n$ means that $V(C_n)=\{v_1,\dots,v_n\}$ and
\[
E(C_n)=\{[v_i,v_{i+1}],i=1,\dots,n-1\}\cup\{[v_1,v_n]\}
\]

Let $G_1$ and $G_2$ be two graphs without common vertices.
We denote by $G_1+G_2$ the graph $G$ for which $V(G)=V(G_1)\cup V(G_2)$ and
$E(G)=E(G_1)\cup E(G_2)\cup E'$ where $E'=\{[x,y]:x\in V(G_1), y\in V(G_2)\}$.

Let $G$ and $H$ be two graphs. We shall say that $H$ is a subgraph of $G$ and
we shall denote $H\se G$ when $V(H)\se V(G)$ and
$E(H)\se E(G)$.

\begin{definition}\label{NN:d11}
A 2-coloring
\beq{NN:eq11}
E(G)=E_1\cup E_2,
\qquad
E_1\cap E_2=\emptyset,
\eeq
is called a blue-red coloring of the edges of the graph $G$ (the edges in $E_1$
are blue and the edges in $E_2$ are red).
\end{definition}

We define for blue-red coloring~\eqref{NN:eq11} and for an arbitrary vertex
$v\in V(G)$
\begin{align*}
N_i(v)&=\{x\in N(v)\mid [v,x]\in E_i\},
\quad
i=1,2;\\
G_i(v)&=G[N_i(v)].
\end{align*}

\begin{definition}\label{NN:d12}
Let $H$ be a subgraph of $G$. We say that $H$ is a monochromatic subgraph
in the blue-red coloring~\eqref{NN:eq11} if $E(H)\se E_1$ or $E(H)\se E_2$.
If $E(H)\se E_1$ we say that $H$ is a blue subgraph and if $E(H)\se E_2$
we say that $H$ is a red subgraph.
\end{definition}

\begin{definition}\label{NN:d13}
The blue-red coloring~\eqref{NN:eq11} is called $(p,q)$-free, if there are
no blue $p$-cliques and no red $q$-cliques. The symbol $G\to(p,q)$ means
that any blue-red coloring of $E(G)$ is not $(p,q)$-free. If $G\to(p,q)$
then $G$ is called edge Folkman $(p,q)$-graph.
\end{definition}

Let $p$, $q$ and $r$ be positive integers. The Folkman number $F(p,q;r)$
is defined by the equality
\[
F(p,q;r)=\min\{|V(G)|:G\to(p,q)\text{ and }\cl<r\}.
\]

In \cite{NN:1} Folkman proved that
\[
F(p,q;r)\text{ exists}\iff r>\max\{p,q\}.
\]

That is why the numbers $F(p,q;r)$ are called Folkman numbers. Only few
Folkman numbers are known. An exposition of the results on the Folkman
numbers was given in \cite{NN:6}. In \cite{NN:6} we computed a new Folkman
number, namely $F(3,4;8)=16$. This result is based upon the fact that
$K_1+C_5+C_5+C_5\to(3,4)$, which was announced without proof in \cite{NN:6}.
In this paper we give a detailed proof of this fact. So, the aim of this paper
is to prove the following

\begin{MT}
Let $G=K_1+\Cf[1]+\Cf[2]+\Cf[3]$, where $\Cf[1]$, $\Cf[2]$,
$\Cf[3]$ are copies of the 5-cycle $C_5$. Then $G\to(3,4)$.
\end{MT}

\section{Auxiliary results}

\begin{lem}\label{NN:l21}
Let $E(G)=E_1\cup E_2$ be a $(3,4)$-free red-blue coloring of the edges of
the graph $G$. Then:
\begin{enumerate}[\rm(a)]
\item
$G_1(v)$ is a red subgraph, $v\in V(G)$;
\item
$(E(G_2(v))\cap E_1)\cup(E(G_2(v))\cap E_2)$ is a $(3,3)$-free red-blue coloring of
$E(G_2(v))$, $v\in V(G)$. Thus $G_2(v)\not\to(3,3)$.
\end{enumerate}
\end{lem}

\begin{proof}
The statement of (a) is obvious. Assume that (b) is not true. Than, since there
is no blue 3-clique, $G_2(v)$ contains a red 3-clique. This red 3-clique
together with the vertex $v$ form a red 4-clique, which is a contradiction.
\end{proof}

\begin{cor}\label{NN:cor21}
Let $E(G)=E_1\cup E_2$ be a $(3,4)$-free blue-red coloring of $E(G)$. Then:
\begin{enumerate}[\rm(a)]
\item
$\cl[G_1(v)]\le 3$,
$v\in V(G)$;
\item
$\cl[G_2(v)]\le 5$,
$v\in V(G)$;
\item
$G_2(v)\nsupseteqq K_3+C_5$, $v\in V(G)$.
\end{enumerate}
\end{cor}

\begin{proof}
The statement of (a) follows from Lemma~\ref{NN:l21}(a). The statements of
(b) and (c) follow from Lemma~\ref{NN:l21}(b), since $K_6\to(3,3)$,
\cite{NN:4} and $K_3+C_5\to(3,3)$, \cite{NN:2}.
\end{proof}

\begin{lem}[\cite{NN:5}]\label{NN:l22}
Let $G=C_5+H$, where $V(H)=\{x,y,z\}$ and $E(H)=\{[x,y],[x,z]\}$. Let
$E(G)=E_1\cup E_2$ be a $(3,3)$-free blue-red coloring of $E(G)$.
Then $H$ is monochromatic in this coloring.
\end{lem}

\begin{lem}[\cite{NN:3}]\label{NN:l23}
Let $G=C_5+K_2$ and $E(G)=E_1\cup E_2$ be a $(3,3)$-free blue-red coloring of
$E(G)$ such that $E(C_5)\se E_i$. Then $E(K_2)\in E_i$.
\end{lem}

\begin{lem}\label{NN:l24}
Let $G=K_1+\Cf[1]+\Cf[2]+\Cf[3]$, where \Cf[1], \Cf[2], \Cf[3] are copies of
the 5-cycle $C_5$ and $V(K_1)=\{a\}$. Let $E(G)=E_1\cup E_2$ be a blue-red
coloring of $E(G)$ such that $\cl[G_1(a)]\le 3$ and $G_2(a)\not\to(3,3)$.
Then, up to numeration of the 5-cycles \Cf[1], \Cf[2] and \Cf[3] we have:
\begin{enumerate}[\rm(a)]
\item
$N_1(a)\supset\VC[1]$ and $N_1(a)\cap\VC[2]$ is an independent set;
\item
$N_2(a)\supset\VC[3]$ and $N_2(a)\cap\VC[2]$ is not an independent set.
\end{enumerate}
\end{lem}

\begin{proof}
Let $\Cf[1]=v_1v_2v_3v_4v_5$, $\Cf[2]=u_1u_2u_3u_4u_5$ and
$\Cf[3]=w_1w_2w_3w_4w_5$. We shall use the following obvious fact
\beq{NN:eq21}
\chi(C_5)=3.
\eeq
It follows from~\eqref{NN:eq21} that
\beq{NN:eq22}
\text{$N_1(a)\cap\VC$ or $N_2(a)\cap\VC$ is not an independent set, $i=1,2,3$.}
\eeq
By~\eqref{NN:eq22} and Corollary~\ref{NN:cor21}(b), at least one of the sets
$N_2(a)\cap\VC$, $i=1,2,3$, is an independent set. Thus, at least one
of the sets $N_1(a)\cap\VC$, $i=1,2,3$, is not an independent set.
Without loss of generality we can assume that
\beq{NN:eq23}
\text{$N_1(a)\cap\VC[1]$ is not an independent set.}
\eeq
It follows from Corollary~\ref{NN:cor21}(a) and \eqref{NN:eq23} that
$N_1(a)\cap\VC[2]=\emptyset$ or $N_1(a)\cap\VC[3]=\emptyset$.
Let for example $N_1(a)\cap\VC[3]=\emptyset$. Then
\beq{NN:eq24}
N_2(a)\supset\VC[3].
\eeq
We have from~\eqref{NN:eq23} and Corollary~\ref{NN:cor21}(a) that
$N_1(a)\cap\VC[2]$ is an independent set. Thus, it follows from~\eqref{NN:eq21}
that $N_2(a)\cap\VC[2]$ is not an independent set. This fact together with
\eqref{NN:eq24} and Corollary~\ref{NN:cor21}(c) give us that
$N_2(a)\cap\VC[1]=\emptyset$. Hence, $N_1(a)\supseteq\VC[1]$.
The Lemma is proved.
\end{proof}

\begin{lem}\label{NN:l25}
Let $G=K_1+\Cf[1]+\Cf[2]+\Cf[3]$, where \Cf, $i=1,2,3$, are
copies of the 5-cycle $C_5$. Let $E(G)=E_1\cup E_2$ be a blue-red coloring
such that some of the cycles \Cf[1], \Cf[2], \Cf[3] is not monochromatic.
Then this coloring is not $(3,4)$-free.
\end{lem}

\begin{proof}
Let $V(K_1)=\{a\}$, $\Cf[1]=v_1v_2v_3v_4v_5$, $\Cf[2]=u_1u_2u_3u_4u_5$ and
$\Cf[3]=w_1w_2w_3w_4w_5$. Assume the opposite, i.e.\ $E(G)=E_1\cup E_2$
is $(3,4)$-free. Then by Corollary~\ref{NN:cor21}(a) we have $\cl[G_1(a)]\le 3$
and by Lemma~\ref{NN:l21}(b) we have $G_2(a)\not\to(3,3)$. Thus, according to
Lemma~\ref{NN:l24} we can assume that
\begin{align}\label{NN:eq25}
N_1(a)&\supseteq\VC[1]\text{ and }N_1(a)\cap\VC[2]
\text{ is independent;}\\
N_2(a)&\supseteq\VC[3]\text{ and }N_2(a)\cap\VC[2]
\text{ is not independent.}
\label{NN:eq26}
\end{align}
It follows from~\eqref{NN:eq25} and Lemma~\ref{NN:l21}(a) that
\beq{NN:eq27}
E(\Cf[1])\se E_2.
\eeq

We have from the statement of the Lemma~2.5 that at least one of the cycles
\Cf, $i=1,2,3$, is not monochromatic and since $E(\Cf[1])\se E_2$
it remains to consider the following two cases:

\emph{Case 1.} \Cf[2] is not monochromatic.
Let for example $[u_1,u_5]\in E_1$ and $[u_1,u_2]\in E_2$.
If $u_1,u_2,u_5\in N_2(a)$ by~\eqref{NN:eq26} we have $G_2(a)\supset C_2^{(3)}+
G[u_1,u_2,u_5]$. It follows from Lemma~\ref{NN:eq22} that $G_2(a)$ contains
a monochromatic 3-clique. This contradicts Lemma~\ref{NN:l21}(b). So, at least
one of the vertices $u_1$, $u_2$, $u_5$ belongs to $N_1(a)$. Therefore, we have
the following subcases:

\emph{Subcase 1a.} $u_1\in N_1(a)$.
Since there are no blue 3-cliques it follows from~\eqref{NN:eq25} that
\beq{NN:eq28}
N_2(u_1)\supset\VC[1].
\eeq
As $[u_1,a]$, $[u_1,u_5]\in E_1$ and $\cl[G_1(u_1)]\le 3$
(see Corollary~\ref{NN:cor21}(a)), the set $N_1(u_1)\cap\VC[3]$
is independent. Therefore, $N_2(u_1)\cap\VC[3]$ is not independent.
This fact together with $[u_1,u_2]\in E_2$ and~\eqref{NN:eq28} give us
$G_2(u_1)\supset K_3+\Cf[1]$, which contradicts Corollary~\ref{NN:cor21}(c).

\emph{Subcase 1b.} $u_2\in N_1(a)$ and $u_1\in N_2(a)$.
Since there are no blue 3-cliques it follows from~\eqref{NN:eq25} that
\beq{NN:eq29}
N_2(u_2)\supset\VC[1].
\eeq
If $N_2(u_1)\cap\VC[1]$ contains two adjacent vertices then these vertices
together with $u_1$ and $u_2$ form a red 4-clique according to~\eqref{NN:eq27}
and~\eqref{NN:eq29}. Hence, $N_2(u_1)\cap\VC[1]$ is independent and, therefore,
$N_1(u_1)\cap\VC[1]$ is not independent. Since $u_5\in N_1(u_1)$ and
$\cl[G_1(u_1)]\le 3$ (see Corollary~\ref{NN:cor21}(a)) we have
$N_1(u_1)\cap\VC[3]=\emptyset$. Hence
\beq{NN:eq210}
N_2(u_1)\supset\VC[3].
\eeq
By~\eqref{NN:eq26} and~\eqref{NN:eq210}
\[
\VC[3]\se N_2(u_1)\cap N_2(a).
\]
Since $[a,u_1]\in E_2$ and there are no red 4-cliques we obtain that
\beq{NN:eq211}
E(\Cf[3])\se E_1.
\eeq
As there are no blue 3-cliques from~\eqref{NN:eq211} it follows that
$N_1(u_2)\cap\VC[3]$ is independent. Therefore, $N_2(u_1)\cap\VC[3]$
contains two adjacent vertices. This fact together with $[u_1,u_2]\in E_2$
and~\eqref{NN:eq29} give us $G_2(u_2)\supset K_3+\Cf[1]$, which contradicts
Corollary~\ref{NN:cor21}(c).

\emph{Subcase 1c.} $u_5\in N_1(a)$ and $u_1,u_2\in N_2(a)$.
Since $a,u_2\in N_2(u_1)$, it follows from Corollary~\ref{NN:cor21}(b) that
at least one of the sets $N_2(u_1)\cap\VC[1]$ and $N_2(u_1)\cap\VC[3]$
is independent. Hence at least one of the sets $N_1(u_1)\cap\VC[1]$,
$N_1(u_1)\cap\VC[3]$, is not independent. Assume that $N_1(u_1)\cap\VC[1]$
is not independent. This fact together with $u_5\in N_1(u_1)$ and
Corollary~\ref{NN:cor21}(a) imply
\beq{NN:eq212}
N_2(u_1)\supset\VC[3].
\eeq
As $[a,u_1,u_2]$ is a red 3-clique and $[a,u_1,u_2,w_i]$ is not a red 4-clique,
$i=1,\dots,5$, it follows from~\eqref{NN:eq26} and~\eqref{NN:eq212} that
$[u_2,w_i]\in E_1$, $i=1,\dots,5$, i.e.\ $N_1(u_2)\supset\VC[3]$. We have from
Lemma~\ref{NN:l21}(a) that $E(\Cf[3])\se E_2$. Thus, according
to~\eqref{NN:eq26} and~\eqref{NN:eq212}, the vertices $a$ and $u_1$ together
with two adjacent vertices of \Cf[3] form a red 4-clique, which is
a contradiction.

Let us now consider the situation when $N_1(u_1)\cap\VC[3]$ is not independent.
Corollary~\ref{NN:cor21}(a) and $u_5\in N_1(u_1)$ imply
\beq{NN:eq213}
N_2(u_1)\supset\VC[1].
\eeq
If $N_2(u_1)\cap\VC[3]\ne\emptyset$ then from $a,u_2\in N_2(u_1)$
and~\eqref{NN:eq213} it follows that $G_2(u_1)\supset K_3+\Cf[1]$,
which contradicts the Corollary~\ref{NN:cor21}(c).
Hence $N_2(u_1)\cap\VC[3]=\emptyset$, i.e.
\beq{NN:eq214}
N_1(u_1)\supset\Cf[3].
\eeq
Since there are no blue 3-cliques we obtain from~\eqref{NN:eq214} and
Lemma~\ref{NN:l21}(a) that
\beq{NN:eq215}
E(\Cf[3])\se E_2.
\eeq
If $N_2(u_2)\cap\VC[3]$ is not independent then according to~\eqref{NN:eq26}
and~\eqref{NN:eq215} an edge in $N_2(u_2)\cap\VC[3]$ together with $a$ and
$u_2$ form a red 4-clique. Let $N_2(u_2)\cap\VC[3]$ be independent. Then
$N_1(u_2)\cap\VC[3]$ is not independent. Thus, it follows from
Corollary~\ref{NN:l21}(a) that $N_1(u_2)\cap\VC[1]$ is independent and
$N_2(u_2)\cap\VC[1]$ is not independent. Then an edge in $N_2(u_2)\cap\VC[1]$
together with the vertices $u_1$ and $u_2$ form a red 4-clique, according
to~\eqref{NN:eq27} and~\eqref{NN:eq213}, which is a contradiction.

\emph{Case~2.} \Cf[3] is not monochromatic but \Cf[2] is monochromatic.
Without loss of generality we can assume that $[w_1,w_5]\in E_1$ and
$[w_1,w_2]\in E_2$. Since $a,w_2\in N_2(w_1)$ it follows from
Corollary~\ref{NN:cor21}(b) that at least one of the sets $N_2(w_1)\cap\VC[1]$
and $N_2(w_1)\cap\VC[2]$ is independent. Hence at least one of the sets
$N_1(w_1)\cap\VC[1]$, $N_1(w_1)\cap\VC[2]$ is not independent.
We shall consider these possibilities:

\emph{Subcase 2a.} $N_1(w_1)\cap\VC[1]$ is not independent.
Since $[w_1,w_5]\in E_1$ it follows from Corollary~\ref{NN:cor21}(a) that
$N_1(w_1)\cap\VC[2]=\emptyset$, i.e.
\beq{NN:eq216}
N_2(w_1)\supset\VC[2].
\eeq
By Lemma~\ref{NN:l21}(b) $G_2(w_1)$ does not contain a monochromatic 3-clique
and $G_2(w_1)\supset\Cf[2]+[a,w_2]$. Since \Cf[2] is monochromatic and
$[a,w_2]\in E_2$, it follows from Lemma~\ref{NN:l23} that
\beq{NN:eq217}
E(\Cf[2])\se E_2.
\eeq
We see from~\eqref{NN:eq26}, \eqref{NN:eq216} and~\eqref{NN:eq217} that
the vertices $a$ and $w_1$ together with an edge of \Cf[2] form a red 4-clique
which is a contradiction.

\emph{Subcase 2b.} $N_1(w_1)\cap\VC[2]$ is not independent.
Since $w_5\in N_1(w_1)$ it follows from Corollary~\ref{NN:cor21}(a) that
\beq{NN:eq218}
N_2(w_1)\supset\VC[1].
\eeq
Corollary~\ref{NN:cor21}(c) and $G_2(w_1)\supset\Cf[1]+[a,w_2]=K_2+\Cf[1]$
imply
\beq{NN:eq219}
N_1(w_1)\supset\VC[2].
\eeq
Lemma~\ref{NN:l21}(a) and~\eqref{NN:eq219} give
\beq{NN:eq220}
E(\Cf[2])\se E_2.
\eeq
Since there are no blue 3-cliques and $[w_1,w_5]\in E_1$ it follows
from~\eqref{NN:eq219} that
\beq{NN:eq221}
N_2(w_5)\supset\VC[2].
\eeq
We see from~\eqref{NN:eq26}, \eqref{NN:eq220} and~~\eqref{NN:eq221} that
the vertices $a$ and $w_5$ together with an edge of \Cf[2] form a red 4-clique
which is a contradiction.
\end{proof}

\section{A property of the graph $C_5+C_5+C_5$}

Let $G=\Cf[1]+\Cf[2]+\Cf[3]$ where \Cf, $i=1,2,3$, are copies of the 5-cycle
$C_5$. Let us consider the blue-red coloring where
$E_1=E(\Cf[1])\cup E(\Cf[2])\cup E(\Cf[3])$. It is clear that this coloring is
$(3,4)$-free. Thus $G\not\to(3,4)$. However the following theorem holds:

\begin{theorem}\label{NN:th31}
Let $G=\Cf[1]+\Cf[2]+\Cf[3]$ where \Cf, $i=1,2,3$, are copies of the $5$-cycle
$C_5$. Let $E(G)=E_1\cup E_2$ be a blue-red coloring such that
$E(\Cf[1])\se E_2$, $E(\Cf[2])\se E_1$ and $E(\Cf[3])\se E_1$.
Then this coloring is not $(3,4)$-free.
\end{theorem}

\begin{proof}
Assume the opposite, i.e.\ that there are no blue 3-cliques and no red 4-cliques.
Let $\Cf[1]=v_1v_2v_3v_4v_5$, $\Cf[2]=u_1u_2u_3u_4u_5$, $\Cf[3]=w_1w_2w_3w_4w_5$.
Since the cycles \Cf[2] and \Cf[3] are blue and there are no blue 3-cliques
we have that the sets $N_1(v_i)\cap\VC[2]$ and $N_1(v_i)\cap\VC[3]$
are independent. Thus, we have
\beq{NN:eq31}
|N_2(v_i)\cap\VC[2]|\ge 3,
\quad
|N_2(v_i)\cap\VC[3]|\ge 3,
\quad
i=1,\dots,5.
\eeq
It follows from~\eqref{NN:eq31} that
\beq{NN:eq32}
N_2(x)\cap N_2(y)\cap\VC\ne\emptyset,
\quad
i=2,3,
\quad
x,y\in\VC[1].
\eeq

Let $x,y\in\VC[1]$. We define
\begin{align*}
B_1(x,y)&=\{v\in\VC[2]\mid[x,v],[y,v]\in E_2\},\\
B_2(x,y)&=\{v\in\VC[3]\mid[x,v],[y,v]\in E_2\}.
\end{align*}
We see from~\eqref{NN:eq32} that
\beq{NN:eq33}
B_i(x,y)\ne\emptyset,
\quad
i=1,2,
\quad
x,y\in\VC[1].
\eeq

We shall prove that
\beq{NN:eq34}
\text{If $[x,y]\in E(\Cf[1])$ then $B_i(x,y)$ is independent, $i=1,2$.}
\eeq

Assume the opposite and let for example $u',u''\in B_1(x,y)$ and
$[u',u'']\in E(\Cf[2])$. By~\eqref{NN:eq33} there exists $w\in B_2(x,y)$.
Since there are no blue 3-cliques then at least one of the edges $[u',w]$,
$[u'',w]$ is red. Hence $[x,y,u',w]$ or $[x,y,u'',w]$ is a red 4-clique,
which is a contradiction.

Let $u'$ and $u''$ be adjacent vertices in \Cf[2]. Since $[u',u'']\in E_1$
and there are no blue 3-cliques we have
\[
N_1(u')\cap N_1(u'')\cap\VC[1]=\emptyset.
\]
Thus $|N_1(u')\cap\VC[1]|\le 2$ or $|N_1(u'')\cap\VC[1]|\le 2$. Hence
\beq{NN:eq35}
|N_2(u')\cap\VC[1]|\ge 3
\text{ and }
|N_2(u'')\cap\VC[1]|\ge 3.
\eeq
So, \eqref{NN:eq35} holds for every two adjacent vertices in \Cf[2].
Hence $|N_2(u)\cap\VC[1]|\ge 3$ holds for at least three vertices in \Cf[2].
Thus, there exist two adjacent vertices in \Cf[2], for example $u_1$ and $u_2$,
such that
\beq{NN:eq36}
|N_2(u_1)\cap\VC[1]|\ge 3
\text{ or }
|N_2(u_2)\cap\VC[1]|\ge 3.
\eeq
If the both inequalities in~\eqref{NN:eq36} are strict then $N_2(u_1)\cap
N_2(u_2)\cap\VC[1]$ contains two adjacent vertices $v'$ and $v''$. Since
$u_1,u_2\in B(v',v'')$ then this contradicts~\eqref{NN:eq34}. Thus,
we may assume that $|N_2(u_1)\cap\VC[1]|=3$. Hence $N_2(u_1)\cap\VC[1]$
contains two adjacent vertices, for example $v_3$ and $v_4$. Now we shall
prove that the third vertex in $N_2(u_1)\cap\VC[1]$ is the vertex $v_1$.
Assume the opposite. Then $v_2\in N_2(u_1)\cap\VC[1]$ or $v_5\in N_2(u_1)\cap
\VC[1]$. Let $v_2\in N_2(u_1)\cap\VC[1]$. Then $v_1,v_5\in N_1(u_1)$. Since
$v_1,v_5,u_2\in N_1(u_1)$ it follows from Corollary~\ref{NN:cor21}(a) that
$N_1(u_1)\cap\VC[3]=\emptyset$. Thus, $G_2(u_1)$ contains $\Cf[3]+[v_3,v_4]K_2+C_5$. According to Lemma~\ref{NN:l21}(b) $G_2(u_1)$
does not contain monochromatic 3-cliques. As $E(\Cf[3])\se E_1$ and
$[v_3,v_4]\in E_2$, this contradicts Lemma~\ref{NN:l23}. We proved that
$v_2\notin N_2(u_1)$. Analogously we prove that $v_5\notin N_2(u_1)$. So,
\beq{NN:eq37}
v_1,v_3,v_4\in N_2(u_1)\text{ and }v_2,v_5\in N_1(u_1).
\eeq
By~\eqref{NN:eq33} we can assume that $w_1\in B_2(v_3,v_4)$. Since
$[v_3,v_4,u_1,w_1]$ is not a red 4-clique we have
\beq{NN:eq38}
[u_1,w_1]\in E_1.
\eeq
As there are no blue 3-cliques and $[u_1,v_2]$, $[u_1,v_5]\in E_1$, it follows
that $[w_1,v_2]$, $[w_1,v_5]\in E_2$. Taking into consideration
$w_1\in B_2(v_3,v_4)$ we have
\beq{NN:eq39}
[w_1,v_i]\in E_2,
\quad
i=2,3,4,5.
\eeq
By~\eqref{NN:eq33} there is $u\in B_1(v_2,v_3)$. Since $[v_2,u_1]\in E_1$
then $u\ne u_1$. We shall prove that $u=u_3$ or $u=u_4$. Assume the opposite.
Then $u=u_2$ or $u=u_5$. Let, for example, $u=u_2$. Since $[v_2,v_3,u_2,w_1]$
is not a red 4-clique, it follows from~\eqref{NN:eq39} and $u_2\in B_1(v_2,v_3)$
that $[u_2,w_1]\in E_1$. We obtained the blue 3-clique $[u_1,u_2,w_1]$ which
is a contradiction. This contradiction proves that $u=u_3$ or $u=u_4$. We can
assume without loss of generality that $u=u_3$. We have
\beq{NN:eq310}
[u_3,w_1]\in E_1,
\eeq
because $[v_2,v_3,u_3,w_1]$ is not a red 4-clique. By~\eqref{NN:eq33} there
exists $u\in B_1(v_4,v_5)$. Repeating the above considerations about
$u\in B_1(v_2,v_3)$ we see that $u=u_3$ or $u=u_4$.

\emph{Case~1.} $u=u_4$.
Since $[v_4,v_5,w_1,u_4]$ is not a red 4-clique, we have $[u_4,w_1]\in E_1$.
Hence $[u_3,u_4,w_1]$ is a blue 3-clique, which is a contradiction.

\emph{Case~2.} $u=u_3$.
In this case we have $u_3\in B_1(v_2,v_3)\cap B_1(v_4,v_5)$, i.e.
\beq{NN:eq311}
[u_3,v_i]\in E_2,
\quad
i=2,3,4,5.
\eeq
As $[v_1,w_1,u_3]$ is not a blue 3-clique, it follows from~\eqref{NN:eq310}
that $[v_1,u_3]\in E_2$ or $[v_1,w_1]\in E_2$.

\emph{Subcase 2a.} $[v_1,u_3]\in E_2$.
By~\eqref{NN:eq311} $N_2(u_3)\supset\Cf[1]$. Since there are no blue 3-cliques
$N_2(u_3)$ contains two adjacent vertices $w',w''\in\VC[3]$. Thus
$G_2(u_3)\supset\Cf[1]+[w',w'']$. By Lemma~\ref{NN:l21}(b) $G_2(u_3)$
contains no monochromatic 3-cliques. This contradicts Lemma~\ref{NN:l23}
because $E(\Cf[1])\se E_2$ and $[w',w'']\in E_1$.

\emph{Subcase 2b.} $[v_1,w_1]\in E_2$.
By~\eqref{NN:eq39} we see that $N_2(w_1)\supset\VC[1]$. Since there are no
blue 3-cliques $N_2(w_1)$ contains two adjacent vertices $u',u''\in\VC[2]$.
Hence $N_2(w_1)\supset\Cf[1]+[u',u'']$ which contradicts Lemma~\ref{NN:l23}.

The theorem is proved.
\end{proof}

\section{Proof of Main Theorem}

Let $\Cf[1]=v_1v_2v_3v_4v_5$, $\Cf[2]=u_1u_2u_3u_4u_5$, $\Cf[3]=w_1w_2w_3w_4w_5$
and $V(K_1)=\{a\}$. Assume the opposite, i.e.\ there exists a $(3,4)$-free
blue-red coloring $E_1\cup E_2$ of the edges of $K_1+\Cf[1]+\Cf[2]+\Cf[3]$.
By Lemma~\ref{NN:l24} we can assume that:
\begin{gather}\label{NN:eq41}
\text{$N_1(a)\supset\VC[1]$ and $N_1(a)\cap\VC[2]$ is independent;}\\
\text{$N_2(a)\supset\VC[3]$ and $N_2(a)\cap\VC[2]$ is not independent.}
\label{NN:eq42}
\end{gather}

We shall prove that
\beq{NN:eq43}
E(\Cf)\se E_2,
\quad
i=1,2,3.
\eeq
By~\eqref{NN:eq41} and Lemma~\ref{NN:l21}(a), $E(\Cf[1])\se E_2$. According
to Lemma~\ref{NN:l25} each of the 5-cycles \Cf[2] and \Cf[3] is monochromatic.
By~\eqref{NN:eq42} $G_2(a)\supset\Cf[3]+e$ where $e\in E(\Cf[2])$.
By Lemma~\ref{NN:l21}(b) $G_2(a)$ contains no monochromatic 3-cliques.
Thus, it follows from Lemma~\ref{NN:l23} that the edge $e$ and the 5-cycle
\Cf[3] have the same color. Therefore, the 5-cycles \Cf[2] and \Cf[3] are
monochromatic of the same color. Thus, it follows from Theorem~\ref{NN:th31}
that $E(\Cf[2])\nsubseteqq E_1$ and $E(\Cf[3])\nsubseteqq E_1$.
We proved~\eqref{NN:eq43}.

Now we shall prove that
\beq{NN:eq44}
N_2(a)=\VC[2]\cup\VC[3].
\eeq
Assume the opposite. Then it follows from~\eqref{NN:eq42} that
$N_1(a)\cap\VC[2]\ne\emptyset$. Let for example $u_1\in N_1(a)\cap\VC[2]$,
i.e.\ $[u_1,a]\in E_1$. We see from~\eqref{NN:eq41} that
\beq{NN:eq45}
[a,u_2]\in E_2.
\eeq
As there are no blue 3-cliques by~\eqref{NN:eq41} and $[u_1,a]\in E_1$
we obtain
\beq{NN:eq46}
N_2(u_1)\supset\VC[1].
\eeq
We see from Corollary~\ref{NN:cor21}(a) that at least one of the sets
$N_2(u_2)\cap\VC[3]$, $N_2(u_2)\cap\VC[1]$ is not independent. If
$N_2(u_2)\cap\VC[1]$ is not independent then it follows from~\eqref{NN:eq46}
and~\eqref{NN:eq43} that the vertices $u_1$ and $u_2$ together with an edge
of \Cf[1] form a red 4-clique. If $N_2(u_2)\cap\VC[3]$ is not independent
then by~\eqref{NN:eq43}, \eqref{NN:eq45} and~\eqref{NN:eq42} the vertices $a$
and $u_2$ together with an edge of \Cf[3] form a red 4-clique.
This contradiction proves~\eqref{NN:eq44}.

It follows from~\eqref{NN:eq44} and Lemma~\ref{NN:l21}(b) that
\beq{NN:eq47}
\text{$\Cf[2]+\Cf[3]$ contains no monochromatic 3-cliques.}
\eeq
Now we obtain from~\eqref{NN:eq47} and~\eqref{NN:eq43}
\begin{gather}\label{NN:eq48}
N_2(x)\cap\VC[3]\text{ is independent},\quad x\in\VC[2];\\
N_2(x)\cap\VC[2]\text{ is independent},\quad x\in\VC[3].
\label{NN:eq49}
\end{gather}

Let us note that
\beq{NN:eq410}
N_1(x)\cap\VC[1]\text{ is independent},\quad x\in\VC[2]\cup\VC[3].
\eeq
Indeed, let for example $x\in\VC[2]$. By~\eqref{NN:eq48} $N_1(x)\cap\VC[3]$
is not independent. This fact and Corollary~\ref{NN:cor21}(a)
prove~\eqref{NN:eq410}.

We shall prove that
\begin{gather}\label{NN:eq411}
N_1(x)\cap\VC[2],\ x\in\VC[1]\text{ is not independent}
\iff
N_2(x)\supset\VC[3];\\
N_1(x)\cap\VC[3],\ x\in\VC[1]\text{ is not independent}
\iff
N_2(x)\supset\VC[2].
\label{NN:eq412}
\end{gather}

The statements~\eqref{NN:eq411} and~\eqref{NN:eq412} are proved analogously.
That is why we shall prove~\eqref{NN:eq411} only. Let $N_1(x)\cap\VC[2]$,
$x\in\VC[1]$ be not independent. Since $[x,a]\in E_1$, it follows from
Corollary~\ref{NN:cor21}(a) that $N_1(x)\cap\VC[3]=\emptyset$, i.e.\
$N_2(x)\supset\VC[3]$. Let now $N_2(x)\supset\VC[3]$, $x\in\VC[1]$. Assume
that $N_1(x)\cap\VC[2]$ is independent. Then $N_2(x)\cap\VC[2]$ is not
independent. Since \Cf[1] is red, $G_2(x)\supset K_3+\Cf[3]$ which contradicts
Corollary~\ref{NN:cor21}(c). So, \eqref{NN:eq411} and~\eqref{NN:eq412}
are proved. Using~\eqref{NN:eq411} and~\eqref{NN:eq412} we shall prove that
\beq{NN:eq413}
\text{$N_1(x)\cap\VC$, $i=2,3$, is independent, $x\in\VC[1]$.}
\eeq
Assume that~\eqref{NN:eq413} is wrong and let for example $N_1(v_1)\cap\VC[2]$
is not independent (remind that $\Cf[1]=v_1v_2v_3v_4v_5$).
Then by~\eqref{NN:eq411} $N_2(v_1)\supset\VC[3]$. If $N_2(v_2)\cap\VC[3]$
is not independent then $v_1$ and $v_2$ together with two adjacent vertices from
$N_2(v_2)\cap\VC[3]$ form a red 4-clique, which is a contradiction. Therefore,
$N_1(v_2)\cap\VC[3]$ is not independent. Thus~\eqref{NN:eq412} gives
$N_2(v_2)\supset\VC[2]$. Repeating the above considerations about the vertex
$v_1$ on $v_2$ we obtain $N_2(v_3)\supset\VC[3]$. In the same way it follows
from $N_2(v_3)\supset\VC[3]$ that $N_2(v_4)\supset\VC[2]$. At the end it follows
from $N_2(v_4)\supset\VC[2]$ that $N_2(v_5)\supset\VC[3]$. So, we proved that
\[
N_2(v_1)\cap N_2(v_5)\supset\VC[3].
\]
Thus, it follows from~\eqref{NN:eq43} that $v_1$ and $v_5$ together with an edge
of \Cf[3] form a red 4-clique, which is a contradiction. This contradiction
proves~\eqref{NN:eq413}. According to~\eqref{NN:eq413} it follows
from~\eqref{NN:eq411} and~\eqref{NN:eq412} that
\beq{NN:eq414}
N_2(x)\not\supset\VC,\ i=2,3,\quad x\in\VC[1].
\eeq
Let $x\in\VC[2]\cup\VC[3]$. By~\eqref{NN:eq410} $|N_1(x)\cap\VC[1]|\le 2$.
Thus, we have the following possibilities:

\emph{Case~1.}
$N_1(x)\cap\VC[1]=\emptyset$ for some vertex $x\in\VC[2]\cup\VC[3]$.
Let for example $N_1(u_1)\cap\VC[1]=\emptyset$ (remind that
$\Cf[2]=u_1u_2u_3u_4u_5$). Then $N_2(u_1)\supset\VC[1]$.
We have from~\eqref{NN:eq410} that $N_2(u_2)\cap\VC[1]$ is not independent.
Thus $u_1$ and $u_2$ together with two adjacent vertices from $N_2(u_2)\cap
\VC[1]$ form a red 4-clique, which is a contradiction.

\emph{Case~2.}
$|N_1(x)\cap\VC[1]|=1$ for some vertex $x\in\VC[2]\cup\VC[3]$.
Let for example $|N_1(u_1)\cap\VC[1]|=1$. Without loss of generality we can
consider that $[u_1,v_1]\in E_1$ and $[u_1,v_i]\in E_2$, $i=2,3,4,5$.
According to~\eqref{NN:eq414} we can assume that $[v_1,w_1]\in E_1$.
Since there are no blue 3-cliques, $[u_1,w_1]\in E_2$. It follows
from~\eqref{NN:eq410} that $N_2(w_1)\cap\VC[1]$ contains two adjacent vertices.
As
\[
N_2(w_1)\cap\VC[1]\se N_2(u_1)\cap\VC[1]=\{v_2,v_3,v_4,v_5\}
\]
we see that $u_1$ and $w_1$ together with two adjacent vertices in
$\{v_2,v_3,v_4,v_5\}$ form a red 4-clique, which is a contradiction.

\emph{Case~3.}
$|N_1(x)\cap\VC[1]|=2$ for every $x\in\VC[2]\cup\VC[3]$.
According to~\eqref{NN:eq48} $N_1(u_1)\cap\VC[3]$ is not independent.
Thus, we can assume that $w_1,w_2\in N_1(u_1)\cap\VC[3]$, i.e.
\beq{NN:eq415}
[u_1,w_1],[u_1,w_2]\in E_1.
\eeq
It follows from~\eqref{NN:eq413}
\beq{NN:eq416}
N_1(w_1)\cap N_1(w_2)\cap\VC[1]=\emptyset.
\eeq
In the considered case we have
\[
|N_1(w_1)\cap\VC[1]|=|N_1(w_2)\cap\VC[1]|=|N_1(u_1)\cap\VC[1]|=2.
\]
We obtain from~\eqref{NN:eq416}
\[
N_1(u_1)\cap N_1(w_1)\cap\VC[1]\ne\emptyset
\text{ or }
N_1(u_1)\cap N_1(w_2)\cap\VC[1]\ne\emptyset.
\]
By~\eqref{NN:eq415} there is a blue 3-clique, which is a contradiction.

Main Theorem is proved.

\section{Example of Folkman edge $(3,5)$-graph without 13-cliques}

Using the Main Theorem we shall prove the following

\begin{theorem}\label{NN:th51}
Let $G=K_4+\Cf[1]+\Cf[2]+\Cf[3]+\Cf[4]$ where \Cf, $i=1,\dots,4$,
are copies of the 5-cycle $C_5$. Then $G\to(3,5)$.
\end{theorem}

In order to prove Theorem~\ref{NN:th51} we shall need the following

\begin{lem}\label{NN:l51}
Let $E(G)=E_1\cup E_2$ is a $(3,5)$-free blue-red coloring of $E(G)$. Then:
\begin{enumerate}[\rm(a)]
\item
$G_1(v)$, $v\in V(G)$, is a red subgraph;
\item
$(E(G_2(v))\cap E_1)\cup(E(G_2(v))\cap E_2)$ is a $(3,4)$-free blue-red coloring
of $E(G_2(v))$, $v\in V(G)$. Thus, $G_2(v)\not\to(3,4)$.
\end{enumerate}
\end{lem}
 Lemma~\ref{NN:l51} is proved in the same way as Lemma~\ref{NN:l21}.

\begin{cor}\label{NN:cor51}
Let $E(G)=E_1\cup E_2$ be a $(3,5)$-free blue-red coloring of $E(G)$. Then:
\begin{enumerate}[\rm(a)]
  \item $\cl[G_1(v)]\le 4$, $v\in V(G)$;
  \item $\cl[G_2(v)]\le 8$, $v\in V(G)$;
  \item $G_2(v)\not\supset K_4+C_5+C_5$;
  \item $G_2(v)\not\supset K_1+C_5+C_5+C_5$.
\end{enumerate}
\end{cor}

\begin{proof}
The statement (a) follows from Lemma~\ref{NN:l51}(a).
The statement (b) follows from Lemma~\ref{NN:l51}(b) and $K_9\to(3,4)$,
\cite{NN:4}. The statement (c) follows from Lemma~\ref{NN:l51}(b) and
$K_4+C_5+C_5\to(3,4)$, \cite{NN:8}. The statement (d) follows from
Lemma~\ref{NN:l51}(b) and Main Theorem.
\end{proof}

\begin{proof}[Proof of Theorem~\ref{NN:th51}]
Assume the opposite, i.e.\ there exists a blue-red coloring $E(G)=E_1\cup E_2$,
which is $(3,5)$-free. Let $V(K_4)=\{a_1,a_2,a_3,a_4\}$.

\emph{Case~1.}
There exists $a_i\in V(K_4)$ such that $|N_1(a_i)\cap V(K_4)|=3$.
Let for example $[a_1,a_2]$, $[a_1,a_3]$, $[a_1,a_4]\in E_1$.
By Corollary~\ref{NN:cor51}(a) at most one of the sets $N_1(a_1)\cap\VC$,
$i=1,2,3,4$, is not empty, i.e.\ $N_2(a_1)$ contains at least three of
the cycles \Cf, $i=1,2,3,4$. Let for example
\[
N_2(a_1)\supset \VC[2]\cup\VC[3]\cup\VC[4].
\]
By Corollary~\ref{NN:cor51}(a) it follows that $N_1(a_1)\cap\VC[1]$
is independent. Thus, $N_2(a_1)\cap\VC[1]\ne\emptyset$. We obtained
that $G_2(a_1)\supset K_1+\Cf[2]+\Cf[3]+\Cf[4]$, which contradicts
Corollary~\ref{NN:cor51}(d).

\emph{Case~2.}
There exists $a_i\in V(K_4)$ such that $|N_1(a_i)\cap V(K_4)|=2$.
Let for example $[a_1,a_2]$, $[a_1,a_3]\in E_1$ and $[a_1,a_4]\in E_2$.
Since $[a_1,a_4]\in E_2$ if the sets $N_2(a_1)\cap\VC$, $i=1,2,3,4$,
are not independent then $G_2(a_1)\supset K_9$, which contradicts
Corollary~\ref{NN:cor51}(b). Hence, at least one of the sets $N_1(a_1)\cap\VC$,
$i=1,2,3,4$, is not independent. Let for example $N_1(a_1)\cap\VC[1]$
is not independent. According to Corollary~\ref{NN:cor51}(a) it follows
from this fact and $[a_1,a_2]$, $[a_1,a_3]\in E_1$ that $N_1(a_1)\cap\VC\emptyset$, $i=2,3,4$, i.e.\ $N_2(a_1)\supset\VC$, $i=2,3,4$. As $[a_1,a_4]\in
E_2$ we have $G_2(a_1)\supset K_1+\Cf[2]+\Cf[3]+\Cf[4]$, which contradicts
Corollary~\ref{NN:cor51}(d).

\emph{Case~3.}
There exist $a_i\in V(K_4)$ such that $|N_1(a_i)\cap V(K_4)|=1$.
Let for example $[a_1,a_2]\in E_1$ and $[a_1,a_3]$, $[a_1,a_4]\in E_2$.
We see from Corollary~\ref{NN:cor51}(a) that at least three of the sets
$N_2(a_1)\cap\VC$, $i=1,2,3,4$, are not independent. Let for example
$N_2(a_1)\cap\VC[2]$, $N_2(a_1)\cap\VC[3]$ and $N_2(a_1)\cap\VC[4]$
are not independent. Since $[a_1,a_3]$, $[a_1,a_4]\in E_2$ it follows
from Corollary~\ref{NN:cor51}(b) that $N_1(a_1)\supset\VC[1]$. According
to Lemma~\ref{NN:l51}(a) it follows from this fact and $[a_1,a_2]\in E_1$
that at least two of the sets $N_1(a_1)\cap\VC$, $i=2,3,4$, are empty.
Therefore, we can assume that $N_2(a_1)\supset\VC[3]$
and $N_2(a_1)\supset\VC[4]$. Since $N_2(a_1)\cap\VC[2]$ is not independent
we have $G_2(a_1)\supset K_4+\VC[3]+\VC[4]$, which contradicts
Corollary~\ref{NN:cor51}(c).

\emph{Case~4.} $E(K_4)\se E_2$.
Since $[a_1,a_i]\in E_2$, $i=2,3,4$, it follows from Corollary~\ref{NN:cor51}(b)
that at least two of the sets $N_1(a_1)\cap\VC$, $i=1,2,3,4$, are not
independent. Let for example $N_1(a_1)\cap\VC[1]$ and $N_1(a_1)\cap\VC[2]$
are not independent. Then by Corollary~\ref{NN:cor51}(a)
$N_1(a_1)\cap\VC[3]=\emptyset$ and $N_1(a_1)\cap\VC[4]=\emptyset$, i.e.
\beq{NN:eq51}
N_2(a_1)\supseteq\VC[3]\cup\VC[4].
\eeq
Since $[a_1,a_i]\in E_2$, $i=2,3,4$, it follows from~\eqref{NN:eq51} and
Corollary~\ref{NN:cor51}(c) that $N_2(a_1)\cap\VC[1]=\emptyset$ and
$N_2(a_1)\cap\VC[2]=\emptyset$. That is why, we have from~\eqref{NN:eq51}
\beq{NN:eq52}
N_1(a_1)=\VC[1]\cup\VC[2].
\eeq
As the vertices $a_1$, $a_2$, $a_3$, $a_4$ are equivalent in this case the above
considerations prove that
\beq{NN:eq53}
\text{$N_1(a_i)$, $i=1,2,3,4$, is a union of two of the cycles \Cf[1], \Cf[2],
\Cf[3], \Cf[4].}
\eeq

Lemma~\ref{NN:l51}(a) and~\eqref{NN:eq52} imply
\beq{NN:eq54}
\Cf[1]+\Cf[2]\text{ is a red subgraph.}
\eeq
Since there are no red 5-cliques we see from~\eqref{NN:eq54} that
\[
N_1(a_i)\cap\VC[1]\ne\emptyset
\text{ or }
N_1(a_i)\cap\VC[2]\ne\emptyset,
\quad
i=2,3,4.
\]
Thus, by~\eqref{NN:eq53} we have that
\beq{NN:eq55}
N_1(a_i)\supset\VC[1]
\text{ or }
N_1(a_i)\supset\VC[2],
\quad
i=2,3,4.
\eeq
Hence, we can assume that
\beq{NN:eq56}
N_1(a_2)\supset\VC[1]
\text{ and }
N_1(a_3)\supset\VC[1].
\eeq

Let $\Cf[1]=v_1v_2v_3v_4v_5$. By~\eqref{NN:eq55} we have the following
possibilities:

\emph{Subcase 4a.} $N_1(a_4)\supset\VC[1]$.
According to~\eqref{NN:eq56} $[v_1,a_i]\in E_1$, $i=1,2,3,4$. Hence,
by Corollary~\ref{NN:cor51}(a) $N_1(v_1)\cap\VC=\emptyset$, $i=2,3,4$, i.e.\
$G_2(v_1)\supset\Cf[2]+\Cf[3]+\Cf[4]$. By~\eqref{NN:eq52} $[v_1,v_2]\in E_2$.
Thus, $G_2(v_1)\supset K_1+\Cf[2]+\Cf[3]+\Cf[4]$, which contradicts
Corollary~\ref{NN:cor51}(d).

\emph{Subcase 4b.}
$N_1(a_4)\cap\VC[1]=\emptyset$, \i.e.\ $N_2(a_4)\supset\VC[1]$.
We have from~\eqref{NN:eq52} and~\eqref{NN:eq56} that $[v_1,a_i]\in E_1$,
$i=1,2,3$, and $[v_1,a_4]\in E_2$. By Corollary~\ref{NN:cor51}(a) at least
two of the sets $N_1(v_1)\cap\VC$, $i=2,3,4$, are empty. Thus, we can assume
that
\beq{NN:eq57}
G_2(v_1)\supset\Cf[3]+\Cf[4].
\eeq
It follows from Corollary~\ref{NN:cor51}(a) that $N_1(v_1)\cap\VC[2]$ is
independent. Hence, $N_2(v_1)\cap\VC[2]$ is not independent. This fact
together with $[v_1,v_2]$, $[v_1,a_4]\in E_2$ and~\eqref{NN:eq57} gives
$G_2(v_1)\supset K_4+\Cf[3]+\Cf[4]$, which contradicts
Corollary~\ref{NN:cor51}(c). This contradiction finishes the proof of
Theorem~\ref{NN:th51}.
\end{proof}

Since $\cl=12$ and $V(G)|=24$ Theorem~\ref{NN:th51} implies

\begin{cor}\label{NN:cor52}
$F(3,5;13)\le 24$.
\end{cor}

Lin proved in~\cite{NN:7} that $F(3,5;13)\ge 18$. In~\cite{NN:9} Nenov
improved this result proving that either $K_8+C_5+C_5\to(3,5)$ or
$F(3,5;13)\ge 19$.

\begin{flushleft}
Faculty of Mathematics and Informatics\\
St.~Kl.~Ohridski University of Sofia\\
5, J.~Bourchier Blvd.\\
BG-1164 Sofia, Bulgaria\\
e-mail: \texttt{nenov@fmi.uni-sofia.bg}
\end{flushleft}

\end{document}